\documentclass[11pt]{amsart}

\title{Orbit inequivalent actions of non-amenable groups}
\author{Inessa Epstein}

\address{Mathematical Sciences Building 6363 \\ University of California \\
    Los Angeles, California~90095 \\ U.S.A.}
\email{iepstein@math.ucla.edu}

\usepackage[margin=1.5in]{geometry}
\usepackage{palatino, amssymb}
\usepackage{verbatim}

\newtheorem{theorem}{Theorem}[section]
\newtheorem{cor}[theorem]{Corollary}
\newtheorem{lemma}[theorem]{Lemma}

\newtheorem{claim}{Claim}[theorem]

\theoremstyle{definition}

\theoremstyle{remark}

\numberwithin{equation}{section}

\newcommand{\N}{{\mathbb N}}
\newcommand{\R}{{\mathbb R}}

\newcommand{\Z}{{\mathbb Z}}
\newcommand{\T}{{\mathbb T}}
\newcommand{\F}{{\mathbb F}}

\newcommand{\mcB}{\mathcal{B}}
\newcommand{\mcC}{\mathcal{C}}

\newcommand{\mcJ}{\mathcal{J}}
\newcommand{\mcI}{\mathcal{I}}

\newcommand{\mcM}{\mathcal{M}}
\newcommand{\mcN}{\mathcal{N}}

\newcommand{\actson}{\curvearrowright}

\newcommand{\ud}{\,\mathrm{d}}

\begin{document}

\begin{abstract} Consider two free measure preserving group actions $\Gamma \actson (X, \mu), \Delta \actson (X,
\mu)$, and a measure preserving action $\Delta \actson^a (Z, \nu)$
where $(X, \mu), (Z, \nu)$ are standard probability spaces. We
show how to construct free measure preserving actions $\Gamma
\actson^c (Y, m)$, $\Delta \actson^d (Y, m)$ on a standard
probability space such that $E_{\Delta}^d \subset E_{\Gamma}^c$
and $d$ has $a$ as a factor. This generalizes the standard notion
of co-induction of actions of groups from actions of subgroups. We
then use this construction to show that if $\Gamma$ is a countable
non-amenable group, then $\Gamma$ admits continuum many orbit
inequivalent free, measure preserving, ergodic actions on a
standard probability space.
\end{abstract}

{\footnote{\emph{Date}: March 2008.
\indent Research partially supported by NSF grant 443948-HJ-21632.
}}

\maketitle

\section{Introduction}

Let us consider a standard probability space $(X, \mu)$ with a
countable group $\Gamma$ acting on $(X, \mu)$ in a Borel measure
preserving manner. This gives rise to the orbit equivalence
relation $E_{\Gamma} = \{ (\gamma \cdot x, x) \mid x \in X \}.$
Two such actions $\Gamma \actson^{a} (X, \mu), \Delta \actson^{b}
(Y, \nu)$ are {\it orbit equivalent} if there exist conull subsets
$A \subset X$, $B \subset Y$ and a measurable, measure preserving
bijection $f \colon A \to B$ such that for any $x, y \in A$, we
have $x E_{\Gamma} y$ if and only if $f(x)E_{\Delta} f(y).$

The theory of orbit equivalence was originally motivated by its
connections to operator algebras. Orbit equivalence first appeared
in a paper by Murray and von Neumann \cite{Murray1936} via the
``group measure space'' construction. One may from a measure
preserving free ergodic action of an infinite countable group
obtain a type II$_1$ von Neumann factor with an abelian Cartan
subalgebra. Two von Neumann algebras obtained in this fashion are
isomorphic via an isomorphism preserving the Cartan subalgebras if
and only if the corresponding actions are orbit equivalent (see
\cite{Feldman1977}).

The first orbit equivalence result is due to Dye \cite{Dye1965},
who showed that all ergodic measure preserving actions of $\Z$ are
orbit equivalent. Later, the work of Ornstein, Weiss, Connes and
Feldman (see \cite{Connes1981}, \cite{Ornstein1980}) provided a
complete classification of ergodic measure preserving actions of
amenable groups. In particular, it was established that all such
actions are orbit equivalent to a $\Z$-action and, consequently,
the orbit equivalence relation remembers only that the group is
amenable.

For non-amenable groups, the situation is quite different. Connes,
Weiss \cite{Connes1980} and Schmidt \cite{Schmidt1981} showed that
all non-amenable groups without Kazhdan's property (T) admit at
least two orbit inequivalent free, measure preserving ergodic
actions. Bezuglyi and Golodets \cite{Bezuglyi1981} showed that
there exists a non-amenable group with continuum many orbit
inequivalent such actions. Results concerning classes of groups
exhibiting this phenomenon of continuum many actions gradually
increased throughout the years. Zimmer \cite{Zimmer1984} showed
that this holds for a number of specific groups with property (T).

Recently, Hjorth \cite{Hjorth2005} showed that actually all groups
with property (T) admit continuum many orbit inequivalent free,
measure preserving, ergodic actions. Gaboriau and Popa
\cite{Gaboriau2005} then used relative property (T) to show this
for all non-cyclic free groups while Ioana \cite{Ioana2006p}
showed this for all groups that admit $\F_2$ as a subgroup. The
question of which groups admit continuum many orbit inequivalent
actions has also been implicitly or explicitly considered as well
as answered for classes of certain groups in the papers
Monod-Shalom \cite{Monod2006}, Popa \cite{Popa2006}, Kechris
\cite{Kechris2005ap}, Tornquist \cite{Tornquist2005}, Fernos \cite
{Fernos2006}.

The main goal of this paper is to present the proof of the following
theorem:

\begin{theorem}\label{main} Let $\Gamma$ be a countable, non-amenable group.
Suppose that there are free, measure preserving actions $\Gamma
\actson (X, \mu)$, $\F_2 \actson (X, \mu)$ on a standard
probability space $(X, \mu)$ such that $\Gamma$ acts ergodically
and $E_{\F_2} \subseteq E_{\Gamma}$. Then $\Gamma$ admits
continuum many orbit inequivalent free, measure preserving,
ergodic actions.
\end{theorem}

Gaboriau and Lyons \cite{Gaboriau2007p} showed that every
countable, non-amenable group admits a free, measure preserving,
ergodic action on a standard probability space $(X, \mu)$ so that
the orbit equivalence relation induced by the action contains the
orbit equivalence relation induced by a free, measure preserving
action of $\F_2$ on $(X, \mu)$. From this and Theorem \ref{main},
we obtain the following corollary:

\begin{cor} Suppose that $\Gamma$ is a countable, non-amenable
group. Then $\Gamma$ induces continuum many orbit inequivalent free,
measure preserving, ergodic actions.
\end{cor}

In \cite{Ioana2006p}, Ioana considered groups $\Gamma$ such that
$\F_2 \leq \Gamma$. Given $\Delta \leq \Gamma$ and an action $a$
of $\Delta$, there is a way to co-induce from this an action of
$\Gamma$ so that the resulting action of $\Gamma$ restricted to
$\Delta$ has the original action by $\Delta$ as a factor. Ioana
then used an action of $\F_2$ on $\T^2$ as well as continuum many
actions of $\F_2$ obtained from irreducible non-isomorphic
representation of $\F_2$ and showed that co-inducing actions of
$\Gamma$ from these actions yields continuum many orbit
inequivalent actions of $\Gamma$. This result uses the fact that
$(\F_2 \ltimes \Z^2, \Z^2)$ has relative property (T) and the fact
that the co-induced action of $\Gamma$ has a strong connection to
the action of $\F_2$ on $\T^2$. Here, the semidirect product $\F_2
\ltimes \Z^2$ is formed by letting SL$_2(\Z)$ act on $\Z^2$ and
viewing $\F_2$ as a finite index subgroup of SL$_2(\Z)$. In
Section 2, we generalize the notion of a co-induction. In
particular, given free measure preserving actions $\Gamma, \Delta
\actson (X, \mu)$ such that $E_{\Delta} \subset E_{\Gamma}$ and a
measure preserving action $\Delta \actson^a (Z, \nu)$, we show how
to construct actions $\Gamma \actson^c (Y, m), \Delta \actson^d
(Y, m)$ so that $E_{\Delta}^d \subset E_{\Gamma}^c$ and $d$ has
$a$ as a factor. In the special case when $\Delta$ is a subgroup
of $\Gamma$, this reduces to the standard induction. In Section
\ref{s:proofthm}, we fit our actions from Section \ref{s:action}
into a theorem of
Ioana and provide the proof of Theorem \ref{main}.\\

Acknowledgements: I would like to thank my advisor, Greg Hjorth,
for his encouragement, support and valuable discussions on the
topic of this paper. I would also like to thank Alexander Kechris
for providing excellent notes and comments concerning the proof I
present here.

\section{The actions of $\Gamma$}\label{s:action}

The aim of this section is to generalize the notion of a
co-induced action.
\\

Throughout this paper, we will use the following notation.

Let $a$ and $b$ be two measure preserving actions of a group
$\Gamma$ on $(X, \mu)$ and $(Y, \nu)$, respectively. $b$ is a {\it
factor} of $a$, written $b \sqsubseteq a$, if there is a Borel
measure-preserving map $p \colon X \to Y$ such that for $\gamma
\in \Gamma$,
    \[p(\gamma^a \cdot x) = \gamma^b \cdot p(x). \]
$E^a_{\Gamma}$ will denote the orbit equivalence relation induced
by $a$ where
    \[E^a_{\Gamma} = \big \{(x, \gamma^{a} \cdot x) \mid x \in X, \gamma \in
    \Gamma \big \}. \]
The {\it diagonal action} $\Gamma \actson^{a \times b}(X \times Y,
\mu \times \nu)$ is given by
    \[\gamma^{a \times b} \cdot (x, y) = (\gamma^{a} \cdot x,
    \gamma^{b} \cdot y). \]
Let $L^2_0(X, \mu) = \{ f \in L^2(X) \mid \int_X f d \mu = 0 \}$.
This is the orthogonal complement of the constant functions in
$L^2(X)$. The {\it Koopman representation} $\kappa_0^a$ of $\Gamma
\actson^{a} (X, \mu)$ on $L^2_0(X)$ is defined by
    \[\gamma^{\kappa_0^{a}}\cdot f(x) = f((\gamma^{-1})^{a} \cdot x). \]
If $\pi_1$ and $\pi_2$ are unitary representations of $\Gamma$, then
$\pi_1 \leq \pi_2$ if $\pi_1$ is isomorphic to a subrepresentation
of $\pi_2$.
    For a Borel space $X$, $\mcB(X)$ will denote the Borel
$\sigma$-algebra on $X$ and $P(X)$ will denote the space of
probability measures on $X$. We will often drop the superscript and
write $\gamma \cdot x$ as opposed to $\gamma^{a} \cdot x$ when it is
clear which action is being used.
\\

First recall the construction of a co-induced action of $\Gamma$
from an action of a subgroup $\Delta$. This first appeared in
\cite{Dooley2008} and can also be found in \cite{Gaboriau2005b}.
Suppose that $\Delta \actson^a (Z, \nu)$ in a Borel measure
preserving manner where $(Z, \nu)$ is a standard probability
space. Let $T \subset \Gamma$ be a left transversal of the cosets
of $\Delta$ in $\Gamma$. Then the space
    \[ Y = \big\{ f \colon \Gamma \to Z \mid f(\gamma_0 \cdot \gamma)
    = \gamma_0^a \cdot f(\gamma) \big\} \]
has a natural identification with the space $Z^{T}$. The
co-induced action of $\Gamma$ on $Z^{T}$ is obtained by
identifying the action of $\Gamma$ on $Y$ given by
  \[ \gamma_0 \cdot f (\gamma) = f(\gamma_0^{-1} \gamma) \]
with the action on $Z^T$ given by
    \[ \gamma \cdot f(t) = \gamma_0^{-1} f(s) \]
where $s \in T$ and $\gamma_0 \in \Delta$ are such that $s
\gamma_0 = \gamma^{-1} t$. This action is then measure preserving
on the standard probability space $(Z^T, \mu^T)$.\\

For our generalization, instead of letting $\Delta \leq \Gamma$,
we assume that the two groups $\Delta$ and $\Gamma$ admit free,
measure preserving actions so that the orbit equivalence relation
of the former is contained in the orbit equivalence relation of
the latter.

\begin{theorem} \label{coinduction} Suppose that there are free measure preserving actions $\Delta \actson^{a_0} (X, \mu)$ and $\Gamma \actson^{b_0}
(X, \mu)$ such that $b_0$ is ergodic and $E_{\Delta}^{a_0} \subset
E_{\Gamma}^{b_0}$. Also, let $\Delta \actson^a (Z, \nu)$ be
measure preserving. Then there is a standard probability space
$(Y, m)$ and actions $c, d$ and map $p$ so that the following
hold:
\begin{enumerate}
    \item  $\Gamma \actson^{c} (Y, m)$ is free, measure
preserving;
    \item $\Delta \actson^{d} (Y, m)$ is free, measure preserving;
    \item $p \colon Y \to Z$ is measure preserving;
    \item \label{l:containment} $E_{\Delta}^d \subset E_{\Gamma}^c;$
    \item \label{l:factormap} $p$ witnesses that $a \sqsubseteq d$.
\end{enumerate}
Moreover, if $a$ is ergodic, then $c$ can be made ergodic as well.
\end{theorem}

\begin{proof} By ergodicity of $\Gamma \actson X$, we may assume that the number
of $\Delta$-equivalence classes in each $\Gamma$-equivalence class
is uniform. In fact, without loss of generality, we will suppose
that each $\Gamma$-equivalence class consists of infinitely many
$\Delta$-equivalence classes.

Consider the space
\[Y = \big\{(x, f) \mid f \colon [x]_{\Gamma} \to Z,
f(\gamma \cdot x_0) = \gamma^{a \times a_{\pi}} \cdot f(x_0) \quad
\forall \gamma \in \Delta \big\},\] which we intend to represent
as the standard probability space $(X\times \Z^{\N}, \mu \times
\nu^{\N})$. In the context of the original co-induced action,
$[x]_{\Gamma}$ takes the place of $\Gamma$. The following lemma is
an adaptation of a coset transversal to our situation.

%EXISTENCE OF TRANSVERSAL
\begin{lemma}There exists a sequence of functions $\{g_i \}_{i \in \N}$ from $X$ to $X$ so that the following
conditions hold:
\begin{enumerate} \label{conditions-on-g}
\item \label{i:gBorel} each $g_i$ is Borel; \item
\label{i:first-element} $g_0(x) = x$ for each $x \in X$; \item
\label{i:transversal} given $x \in X,$ $\{g_i(x) \}_{i \in \N}$
enumerates a transversal for the $\Delta$-equivalence classes in
$[x]_{\Gamma};$ \item \label{i:selector} if $i \neq j$ and $x \in
X$, then $g_i(x) \neq g_j(x)$.
\end{enumerate}
\end{lemma}
\begin{proof} Let $\{ \gamma_n \}_{n
\in \N}$ enumerate the elements of $\Gamma$ so that $\gamma_0 =
e$. We define $g_i \colon X \to X$ inductively on $i$. First,
define
\begin{align*}
h_i(x) &= \hbox{ least } k \in \N \hbox { such that }\\
&\forall l < k, \exists j < i \big((\gamma_l \cdot x,
\gamma_{h_j}(x)) \in
E_{\Delta}^{a_0} \big) \\
&\wedge \quad \forall j < i, \big((\gamma_k \cdot x,
\gamma_{h_j}(x)) \notin E_{\Delta}^{a_0} \big).
\end{align*}
That is, we want to take the least $k$ such that the
$\Delta$-equivalence class of $\gamma_k \cdot x$ has not already
appeared for a previous $h_j$ and the $\Delta$-equivalence class
of each $\gamma_l \cdot x$ for each $l < k$ has already appeared.
Then let
\[g_i(x) = \gamma_{h_i(x)} \cdot x.\]
Conditions \eqref{i:first-element}, \eqref{i:transversal} and
\eqref{i:selector} are clearly satisfied by our construction. It
suffices to check that $h_i$ is Borel since then $g_i$ is a
composition of two Borel maps.

Note that
\begin{align*}
h_i^{-1}(k) =& \big(\bigcup_{l < k} \bigcap_{j < i} \{x \in X \mid
(\gamma_l \cdot x, g_j(x)) \in E_{\Delta}^{a_0} \} \big)\\
&\quad \bigcap (\bigcap_{j < i} \{ x \in X \mid (\gamma_k \cdot x,
g_j(x)) \notin E_{\Delta}^{a_0} \} \big)
\end{align*}
\[\{x \in X \mid (x, \gamma_l \cdot g_j(x)) \in E_{\Delta}^{a_0}\} =
\hbox{proj}_X(E_{\Delta}^{a_0} \cap \{(x, \gamma_l \cdot g_j(x))
\mid x \in X \} )\] and projections of sets with countable
sections are Borel.
\end{proof}

Thus, we have an isomorphism $F \colon Y \to X \times Z^{\N}$
given by
    \[F(x, f) = (x, f(g_0(x)), f(g_1(x)), ...) \]
and we may let $m$ be the product measure $\mu \times \nu^{\N}$.

As for the actions, let $\Gamma \actson^{c} (Y, m)$ be defined by
    \[\gamma^{c} \cdot (x, f) = (\gamma \cdot x, f) \quad \forall \gamma \in
    \Gamma \]
and $\Delta \actson^{d} (Y, m)$ be defined analogously by
    \[ \gamma_0^d \cdot (x, f) = (\gamma_0 \cdot x, f) \quad
    \forall \gamma_0 \in \Delta .\]
\\

We will write the action $c$ as a skew-product action on $X \times
Z^{\N}$ consistent with the above representation. For this
purpose, let $S_{\infty}$ act on $\Delta^{\N}$ by shift, i.e., for
$\alpha \in S_{\infty}$, $\delta \in \Delta^{\N}$
    \[ \alpha \cdot \delta(k) = \delta(\alpha^{-1}  (k)) \]
and consider the semidirect product $S_{\infty} \ltimes
\Delta^{\N}$.
    Then define the cocycles
    \[ \alpha \colon \Gamma \times X \to S_{\infty}, \quad \delta \colon \Gamma \times X \to \Delta^{\N}, \quad \beta \colon \Gamma \times X \to  S_{\infty} \ltimes \Delta^{\N}
\] by
    \[ \alpha(\gamma, x)(k) = n \iff g_k(x) = g_n(\gamma \cdot x) \]
    \[ \delta (\gamma, x)(k)g_{\alpha(\gamma, x)^{-1}(k)}(x) = g_k(\gamma \cdot x) \]
    \[ \beta = (\alpha, \delta) .\]
Then let $S_{\infty} \ltimes \Delta^{\N} \actson Z^{\N}$ by
    \[ (\alpha, \delta) \cdot f(k) = \delta(k) \cdot f(\alpha^{-1} (k)) .\]

\begin{lemma} \label{cocycle} The following hold:
\begin{enumerate}
\item \label{beta-cocycle} $\beta$ is a cocyle; \item
\label{beta-action}$S_{\infty} \ltimes \Delta^{\N} \actson Z^{\N}$
defines an action; \item \label{beta-rep}for all $\gamma \in
\Gamma$ and $(x, f) \in X \times Z^{\N}$,
    \[ \gamma^c \cdot (x, f) = (\gamma \cdot x, \beta(\gamma, x)
\cdot f) .\]
\end{enumerate}

\end{lemma}
\begin{proof}
\eqref{beta-cocycle} Let $\gamma_1, \gamma_2 \in \Gamma$ and $x
\in X$. Observe that $\alpha$ is a cocycle, i.e.,
    \[ \alpha(\gamma_1 \gamma_2, x) = \alpha(\gamma_1, \gamma_2 \cdot
x)\alpha(\gamma_2, x) .\] Indeed, if $\alpha(\gamma_2, x)(k) = l$
and $\alpha(\gamma_1, \gamma_2 \cdot x)(l) = n$, then
    \[ g_k(x)E_{\Delta} g_l(\gamma_2 \cdot x), \quad g_l(\gamma_2 \cdot x) E_{\Delta} g_n(\gamma_1 \gamma_2 \cdot x)
.\] Consequently, $g_k(x)E_{\Delta} g_n(\gamma_1 \gamma_2 \cdot
x)$ and, by the definition of $\alpha$, $\alpha(\gamma_1 \gamma_2,
x)(k) = n$.
\\

It remains to show that
    \[ \beta(\gamma_1 \gamma_2, x) = \beta(\gamma_1, \gamma_2 \cdot
    x)\beta(\gamma_2, x) .\]
By our definitions of $\alpha$ and $\delta$,
\begin{align*}
&\delta(\gamma_1 \gamma_2, x)(k) g_{\alpha(\gamma_1 \gamma_2,
x)^{-1}(k)}(x) \cr \qquad &= g_k(\gamma_1 \gamma_2 \cdot x) \cr
\qquad &= \delta (\gamma_1, \gamma_2 \cdot x)(k) \big
[g_{\alpha(\gamma_1, \gamma_2 \cdot x)^{-1}(k)}(\gamma_2 \cdot x)
\big] \cr \qquad&= \delta (\gamma_1, \gamma_2 \cdot x)(k) \big
[\delta(\gamma_2, x)( \alpha(\gamma_1, \gamma_2 \cdot x)^{-1}(k))
g_{\alpha(\gamma_2, x)^{-1}(\alpha(\gamma_1 ,\gamma_2 \cdot
x)^{-1}(k))}(x) \big] \cr \qquad &= \big[ \delta (\gamma_1,
\gamma_2 \cdot x)(k) \delta(\gamma_2, x)( \alpha(\gamma_1,
\gamma_2 \cdot x)^{-1}(k)) \big] g_{\alpha(\gamma_1 \gamma_2,
x)^{-1}(k)}(x)
\end{align*}
and, as a result,
    \[ \delta(\gamma_1 \gamma_2, x)(k) = \delta (\gamma_1, \gamma_2
\cdot x)(k) \delta(\gamma_2, x)( \alpha(\gamma_1, \gamma_2 \cdot
x)^{-1}(k)).\]

Finally, from the above calculations,
\begin{align*}
&(\alpha(\gamma_1, \gamma_2 \cdot x), \delta(\gamma_1, \gamma_2
x))(\alpha(\gamma_2, x) \delta(\gamma_2, x))\cr \qquad &=
(\alpha(\gamma_1, \gamma_2 \cdot x), \delta(\gamma_1, \gamma_2
x))(\alpha(\gamma_2, x) \delta(\gamma_2, x)) \cr \qquad &=
(\alpha(\gamma_1 \gamma_2, x), \delta(\gamma_1, \gamma_2
x)\big[\alpha(\gamma_1, \gamma_2 x) \cdot \delta(\gamma_2 \cdot x)
\big] \cr \qquad &= (\alpha(\gamma_1 \gamma_2, x), \delta(\gamma_1
\gamma_2, x))
\end{align*}
establishing that $\beta$ is a cocycle.
\\

\eqref{beta-action}
    Let $(\alpha_1, \delta_1), (\alpha_2, \delta_2) \in S_{\infty}
\ltimes \Delta^{\N}$ and $f \in Z^{\N}$. Then
\begin{align*}
(\alpha_1, \delta_1) \cdot (\alpha_2, \delta_2) f(k) &= (\alpha_1,
\delta_1) \cdot \delta_2(k) f(\alpha_2^{-1} (k)) \cr &=
\delta_1(k) \delta_2(\alpha_1^{-1}(k)) f(\alpha_2^{-1}
\alpha_1^{-1}(k)) \cr &= (\alpha_1 \alpha_2, \delta_1(\alpha_1
\cdot \delta_2)) \cdot f(k).
\end{align*}
\\

\eqref{beta-rep} Given $(x, f) \in Y$ and $\gamma \in \Gamma$,
\begin{align*}
f(g_k(\gamma \cdot x)) &= f\big(\delta(\gamma, x)(k)
g_{\alpha(\gamma, x)^{-1}(k)}(x)\big) \cr &= \delta(\gamma, x)(k)
f \big(g_{\alpha(\gamma, x)^{-1}(k)}(x) \big).
\end{align*}

Thus, $\gamma^c(x, f) = (\gamma \cdot x, \beta(\gamma, x) \cdot
f)$.

\end{proof}

Since the actions $\Gamma \actson (X, \mu)$ and $S_{\infty}
\ltimes \Delta^{\N} \actson (Z^{\N}, \nu^{\N})$ are measure
preserving, the action
    \[ \Gamma \actson^c (X \times Z^{\N}, \mu \times \nu^{\N}) \]
formed by a skew-product is measure preserving as well (see
\cite{Glasner2003}).
\\

For the action $\Delta \actson^d (Y, m)$, define the cocycle
$\sigma \colon \Delta \times X \to \Gamma$ by
    \[ \sigma(\gamma_0, x) = \gamma \iff \gamma_0 \cdot x = \gamma
\cdot x \] for any $\gamma_0 \in \Delta$. Then define $\Delta
\actson^d X \times Z^{\N}$ by
    \[ \gamma_0 \cdot (x, f) = \sigma(\gamma_0, x)^c \cdot (x, f) .\]

At this point, we can see that our particular construction
mandates the freeness of the actions $a_0$ and $b_0$ to define the
cocycles $\delta$ and $\sigma$, respectively.

\begin{lemma}\label{c:F2measurepreserve} If $A \subset Y$ is $\Delta$-invariant with $\Gamma$-invariant probability measure $m'$, then the action $\Delta \actson^{d|A}
(A, m')$ is measure preserving.
\end{lemma}
\begin{proof} Let $B \subset A$ be Borel and let $\gamma_0 \in
\Delta$. Then

\begin{align*}
m'(\gamma_0 \cdot B) &= m'\Big(\gamma_0 \cdot \bigcup_{\gamma \in
\Gamma} \big \{ (x, f) \in B \mid \gamma_0 \cdot (x, f) = \gamma
\cdot
(x, f) \big \}\Big) \\
&= m' \Big(\{\bigcup_{\gamma \in \Gamma} \gamma \cdot \big\{ (x,
f) \in
B \mid \gamma_0 \cdot (x, f) = \gamma \cdot (x, f) \big\} \Big)\\
&= \sum_{\gamma \in \Gamma} m' \Big( \big \{ (x, f) \in B \mid
\gamma_0 \cdot (x, f) = \gamma \cdot (x, f) \big \}\Big) = m'(B).
\end{align*}
\end{proof}

Define the map $p \colon Y \to Z$ by $p(x, f) = f(x) = g_0(x)$.
Since for each $(x, f) \in Y$, $f$ is $\Delta$-equivariant, it is
clear that $p$ is also a $\Delta$-equivariant map. To see that $p$
is measure preserving and does, in fact, witness that $a
\sqsubseteq d$, let $A \subset Z$ be arbitrary. Then
    \[ p_*(\mu \times \nu^{\N})(A) = \mu \times \nu^{\N} \Big( \big\{ (x, f) \in Y \mid f(1) \in A
    \big\} \Big )= \nu(A). \]

We will show how to obtain ergodicity of $c$ in the proof of Lemma
\ref{l:action} since we will use some facts concerning ergodic
decomposition we have yet to prove.

\end{proof}

We now need a general lemma concerning ergodic decompositions (see
\cite{Kechris2004}).\\

Let $\Gamma \actson (X, \mu), (Y, \nu)$ be Borel and measure
preserving where $X$ and $Y$ are standard probability spaces and
$\nu$ is ergodic. Suppose that $p \colon X \to Y$ is a
$\Gamma$-equivariant map, i.e.,
    \[ \forall \gamma \in \Gamma \quad p(\gamma \cdot x) = \gamma
\cdot p(x) .\]
Consider the ergodic decomposition of $X$ with
respect to the action $\Gamma \actson X$. This is given by a
$\Gamma$-invariant Borel map $\Phi \colon X \to \mcI$ where $\mcI$
is a standard Borel space and a Borel map $i \in \mcI \mapsto \mu_i
\in P(X)$ such that the following hold:
\begin{enumerate}
    \item for each $i \in \mcI$, if we let
\[X_i = \{x \in X \mid \Phi(x) = i \},\]
then $X_i$ is $\Gamma$-invariant and $\mu_i$ is the unique ergodic
$\Gamma$-invariant measure on $X_i;$
    \item  $\mu = \int_{\mcI} \mu_i \ud \eta(i)$ where $\eta = \Phi_* \mu$.
\end{enumerate}

\begin{lemma} \label{ergdecomp} The following hold:
\begin{enumerate}
\item \label{ed:invariant} If $A \subset X$ is a $\Gamma$-invariant subset and $B \subset
Y$, then
        \[\mu(A \cap p^{-1}(B)) = \nu(B) \mu(A). \]
\item \label{ed:measurepreserve} If $\Delta \actson (X, \mu)$ is
another Borel measure preserving action such that for any
$\Delta$-invariant set $A \subset X$ and any $B \subset Y$, we
have
    \[ \mu(A \cap p^{-1}(B)) = \nu(B) \mu(A) ,\]
then for $\Phi_*\mu$-almost every $i \in \mcI$, if $A \subset X$
is $\Delta$-invariant, we have
    \[ \mu_i(A \cap p^{-1}(B)) = \nu(B) \mu_i(A). \]
\end{enumerate}
\end{lemma}
\begin{proof}

\eqref{ed:invariant} It suffices to show that for some subset
$\mcI_0 \subset \mcI$ such that $\Phi_*\mu(\mcI_0) = 1$, we have
$p_*\mu_i = \nu$ for all $i \in {\mathcal I}_0$. Granted this, we
may finish the proof. Indeed, since $A$ is $\Gamma$-invariant,
then up to null sets, $A = \Phi^{-1}(I)$ for some subset $I
\subset \mcI_0$. Thus,
\begin{align*}
\mu(A \cap p^{-1}(B)) &= \int_I \mu_i(p^{-1}(B)) \ud \Phi_*\mu(i)\\
&=\int_I \nu(B) \ud \Phi_*\mu(i) \\
&=\nu(B) \mu(A).
\end{align*}

Since the measure $\nu$ on $Y$ is ergodic and $\Gamma$-invariant,
we may let $C \subset Y$ be such that $\nu(C) = 1$ and $\nu$ is
the unique $\Gamma$-invariant probability measure on $C$. By the
fact that $p$ is measure-preserving, we have that $\mu(p^{-1}(C))
= 1$. Then for $\Phi_*\mu$-conull many $i \in \mcI$, we have
$p_*\mu_i(C) = 1$. Also, by equivariance of the map $p$,
$p_*\mu_i$ is a $\Gamma$-invariant measure on $C$. By uniqueness
of $\nu$, it must be that $p_*\mu_i = \nu$.

\eqref{ed:measurepreserve} Suppose that this fails on a set of
$\Phi_*\mu$-positive measure. Then, without loss of generality, we
may find a subset $D \subset \mcI$ of $\Phi_*\mu$-positive measure
such that for each $i \in D$, there is a $\Delta$-invariant
$F_{\sigma}$ subset $A_i \subset X$ such that
    \[\mu_i(p^{-1}(B) \cap A_i) < \nu(B) \mu_i(A_i).\]

We show that there is a $\Phi_*\mu$-measurable assignment $\psi$
from $D$ to the $F_{\sigma}$ subsets of $X$ so that for each $i
\in D$, the above inequality holds where $A_i = \psi(i)$.

Let $F \colon \N^{\N} \to F_{\sigma}(X)$ (where $F_{\sigma}(X)$ is
the set of $F_{\sigma}$ subsets of $X$) be a Borel bijection in
the sense that
    \[ \{(w, x) \mid x \in F(w) \} \subset \N^{\N} \times Y \]
is Borel. Then let $F_1 \colon \N^{\N} \to \mcB(X)$ be defined by
$F_1(w) = F(w) \cap p^{-1}(B)$.

Consider the subset $D_0 \subset \N^{\N} \times P(X)$ defined by
\[D_0 = \{(w, \lambda) \in \N^{\N} \times P(X) \mid \lambda(p^{-1}(B) \cap F(w))
< \nu(B)\lambda(F(w)) \}.\] We observe that $D_0$ is Borel. Indeed,
by 17.25 of \cite{Kechris1995}, the maps
    \[(w, \lambda) \in \N^{\N} \times P(X) \mapsto \nu(B)\lambda(F(w)),\]
    \[(w, \lambda) \in \N^{\N} \times P(X) \mapsto \lambda(F_1(w))\]
are Borel and, as a result, the map
    \[(w, \lambda) \in \N^{\N} \times P(X) \mapsto \lambda(p^{-1}(B) \cap F(w)) - \nu(B)\lambda(F(w))\]
is also Borel. This establishes that $D_0$ is Borel. Now let
\[D_1 = \{(w, i) \in \N^{\N} \times \mcI \mid (w, \mu_i) \in D_0 \}.\]
$D_1$ is Borel as well since the assignment $i \in \mcI \mapsto
\mu_i \in P(X)$ is Borel.

Finally, $D = \hbox{proj}_{\mcI}(D_1)$ and $D$ is analytic. Then
$\N^{\N} \times D \subset \N^{\N} \times \mcI$ is analytic as well
and, by 18.1 of \cite{Kechris1995}, there is a
$\Phi_*\mu$-measurable assignment $\psi \colon D \to \N^{\N}$ such
that
    \[(\psi(i), i) \in D_1 \quad \forall i \in D.\]
Note that $\bigcup_{i \in D} A_{\psi(i)}$ is a measurable
$\Gamma$-invariant subset of $X$ of $\mu$-positive measure so we
aim to obtain a contradiction to the fact that
    \[ \mu\big((\bigcup_{i \in D}A_{\psi(i)}) \cap p^{-1}(B)\big) = \nu(B) \mu\big(\bigcup_{i \in
    D}A_{\psi(i)}\big). \]

We have

\begin{align*}
\mu(p^{-1}(B) \cap \bigcup_{i \in D}A_{\psi(i)}) &= \int_D
\mu_i(p^{-1}(B) \cap
A_{\psi(i)}) \ud \Phi_*\mu(i) \\
&< \int_D \nu(B)\mu_i (A_{\psi(i)})  \ud \Phi_*\mu(i) \\
&= \nu(B) \mu\big(\bigcup_{i \in D} A_{\psi(i)} \big).
\end{align*}
\end{proof}

We are now ready to specifically consider an action of $\Gamma$
induced from an action of $\F_2$. Fix actions $\Gamma
\actson^{b_0} (X, \mu)$ and $\F_2 \actson^{a_0} (X, \mu)$
satisfying the hypotheses of Theorem \ref{main}. Let SL$_2(\Z)
\actson (\T^2, h)$ where $h$ is the Haar measure as follows:
    \[ A \cdot (t_1, t_2) = (A^{-1})^t
    \left(\begin{array}{c} t_1 \\ t_2 \end{array}\right) .\]
Fix a copy of $\F_2$ in SL$_2(\Z)$ with finite index and let $a$
be the action of $\F_2$ on $(\T^2, h)$ given by restricting the
action of SL$_2(\Z)$ on $(\T^2, h)$. This action is then free,
measure preserving and weakly  mixing. For more on this, see
Section 16 of \cite{Kechris1995} or \cite{Tornquist2005}.

\begin{lemma} \label{l:action} Given $a$ as specified above, suppose we have the following:
\begin{enumerate}
    \item  $\F_2 \actson^{a_{\pi}}  (Z, \nu)$ is a weakly mixing
action;
    \item $\F^2 \actson^{a \times a_{\pi}} \T^2 \times Z$ is the diagonal
action obtained from $a$ and $a_{\pi}$;
    \item $q \colon \T^2 \times Z \to \T^2$ is given by $q(t, z) = t$.
\end{enumerate}
Then there is a standard probability space $(Y, m)$ and actions
$c, d$ and map $p$ so that the following hold:
\begin{enumerate}
    \item  $\Gamma \actson^{c} (Y, m)$ is free, measure
preserving, ergodic;
    \item $\F_2 \actson^{d} (Y, m)$ is free, measure preserving;
    \item $p \colon Y \to \T^2 \times Z$ is measure preserving;
    \item \label{l:containment} $E_{d} \subset E_{c};$
    \item \label{l:factormap} for any non non-null $d$-invariant
    subset $Y_0 \subset Y$, $p|Y_0$ witnesses that
    \[a \times a_{\pi} \sqsubseteq d|Y_0;\]
\item \label{l:noclasscollapse} $\forall \gamma \in \Gamma \setminus \{ e \}$,
    \[m\Big( \big\{y \in Y \mid q \circ p(\gamma^{c} \cdot y) = q \circ p(y)\big\}\Big) = 0. \]
\end{enumerate}
\end{lemma}

\begin{proof}Consider the space
    \[ (Y, m) = \big (X \times (\T^2 \times Z)^{\N}, \mu \times (h \times \nu)^{\N}\big). \]
We may obtain the actions $c$ and $d$ on $(Y, m)$ from the
construction in Theorem \ref{coinduction}. Note that since $a_0$
and $b_0$ are free, the actions $c$ and $d$ are free as well. We
just need to select a measure on $Y$ so that the action of
$\Gamma$ on $Y$ with respect to this measure is measure preserving
and ergodic. For this purpose, we take an ergodic decomposition of
$Y$ with respect to the action $c$ and let $\Phi \colon Y \to
\mcI$ and $i \in \mcI \mapsto m_i \in P(Y)$ be the corresponding
Borel assignments. Our remaining goal is to show that
$\Phi_*m$-almost every measure in $\mcI$ satisfies our conditions.

\begin{lemma}\label{gamma-decomposition} The following hold:
\begin{enumerate}
\item \label{measurepreserve} If $B \subset \T^2 \times Z$, then
for $\Phi_*m$-almost every $i \in \mcI$, for every
$\F_2$-invariant $A \subset Y$,
     \[ m_i(A \cap p^{-1}(B)) = h \times \nu(B) m_i (A) .\]

\item \label{noncollapse} If $\gamma \in \Gamma \setminus \{ e
\}$, then for almost every $i \in \mcI$,
\[m_i\Big(\big\{(x, f) \in Y_i \mid q \circ p(x, f) = q \circ p(\gamma
\cdot (x, f)) \big\} \Big) = 0. \]
\end{enumerate}
\end{lemma}
\begin{proof}
\eqref{measurepreserve} By Lemma \ref{ergdecomp}, if $A \subset Y$
is $\F_2$-invariant, then
    \[ m(A \cap p^{-1}(B)) = h \times \nu(B) m(A) .\]
Moreover, for $\Phi_*m$-almost every $i \in \mcI$,
    \[ m_i(A \cap p^{-1}(B)) = h \times \nu(B) m_i (A) .\]

\eqref{noncollapse} We first show that for any $\gamma \in \Gamma
\setminus \{ e \}$,
    \[m \Big( \big\{(x, f) \in Y \mid q \circ p(x, f) = q \circ p(\gamma
    \cdot (x, f)) \big\} \Big) = 0. \]

Fix $\gamma \in \Gamma \setminus \{ e \}$. Define
\[Y_{\gamma} = \Big\{(x, f) \in Y \mid q \circ p(x, f) = q \circ p(\gamma
\cdot (x, f)) \Big\} \] and suppose that $m(Y_{\gamma})>0$.
Without loss of generality, we may assume that for all $(x, f) \in
Y_{\gamma}$,
\[g_0(x) = \gamma_0 \cdot g_k(\gamma \cdot x))\]
for some fixed $k \in \N$ and $\gamma_0 \in \F_2$. If $k = 0$, then
$\gamma_0 \neq e$ and since the action of $\F_2$ is free on $\T^2$,
it is impossible that $q(f(x)) = q(f(\gamma \cdot x))$. On the other
hand, if $k \neq 0$, then
\[Y_{\gamma} = \Big \{(x, f) \in Y \mid f(g_0(x)) = \gamma_0 \cdot
f(g_k(x)) \Big \}.\] The measure $h \times \nu$ is non-atomic and,
hence, $m(Y_{\gamma}) = 0.$

Let $D \subset \mcI$ be a set of $\Phi_*m$-positive measure so
that for $i \in D$,
\[m_i\Big( \big\{(x, f) \in Y_i \mid q \circ p(x, f) = q \circ p(\gamma
\cdot (x, f)) \big\} \Big) > \delta \]
for some $\delta > 0$. We've established $Y_{\gamma}$ has measure
zero. Thus,
\begin{align*}
0 &= m(Y_{\gamma} \cap \bigcup_{i \in D} Y_i) \\
&= \int_D m_i(Y_{\gamma} \cap Y_i) \ud \Phi_*m(i) \\
&> \delta \int_D \ud \Phi_*m(i) \\
&= \delta \cdot m(\bigcup_{i \in D} Y_i)
\end{align*}
which is impossible by our choice of $D$.

\end{proof}

Let $\mcB = \{ B_n \}_{n \in \N}$ generate the Borel
$\sigma$-algebra on $\T^2 \times Z$. Without loss of generality,
we may assume that $\mcB$ is clopen, invariant under the action of
$\F_2$ and closed under Boolean operations. Let $i \in \mcI$ be
such that conditions \eqref{measurepreserve} and
\eqref{noncollapse} of Lemma \ref{gamma-decomposition} hold for
all $\gamma \in \Gamma \setminus \{e\}$ and $B \in \mcB$. Since
$Y_i$ is $\Gamma$-invariant, it is $\F_2$-invariant as well and,
as a result, $E_{d|Y_i} \subset E_{c|Y_i}$. $\F_2 \actson^{d|Y_i}
(Y_i, m_i)$ is measure preserving by Lemma
\ref{c:F2measurepreserve}, conditions \eqref{l:factormap} and
\eqref{l:noclasscollapse} follow from Lemma
\ref{gamma-decomposition} \eqref{measurepreserve} and Lemma
\ref{gamma-decomposition} \eqref{noncollapse}, respectively. Thus,
$(Y_i, m_i)$ with actions $\Gamma \actson^{c|Y_i} (Y_i, m_i)$ and
$\F_2 \actson^{d|Y_i} (Y_i, m_i)$ and factor map $p$ are as
desired.

Looking back at Theorem \ref{coinduction}, note that Lemma
\ref{ergdecomp} makes no assumptions on the action of $\Delta$ on
$(Y, \nu)$ except ergodicity and neither does Lemma
\ref{gamma-decomposition} \eqref{measurepreserve} on the action of
$\F_2$ on $(\T^2 \times Z, h \times \nu)$. Thus, the proof here
shows that in Theorem \ref{coinduction}, the action of $c$ can be
made ergodic when $a$ is ergodic.

\end{proof}

\section{Proof of the main theorem}\label{s:proofthm}

We now proceed as in Ioana \cite{Ioana2006p} with the action
described in Section \ref{s:action} replacing the co-induced
action and making a change to the order of operations in
constructing actions.
\\

As defined in the previous section, $a$ is the action of $\F_2$ on
$(\T^2, h)$. The following lemma is Theorem 1.3 of
\cite{Ioana2006p}:

\begin{lemma} Let $\Gamma$ be a group such that $\F_2 \leq \Gamma$ is a subgroup and
let $\{c_i \}_{i \in I}$ be an uncountable set of orbit equivalent
ergodic, free, measure preserving actions $\Gamma \actson^{c_i}(Y_i,
m_i)$ so that the following conditions hold:
\begin{enumerate}
\item \label{cond:quotient} $a$ is a quotient of $c_i|_{\F_2}$ with
quotient map $p_i: Y_i \rightarrow \T^2;$
\item \label{cond:noclasscollapse} $\forall i \in I \quad \forall \gamma \in \Gamma
\setminus \{ e \}$,
\[m\Big( \big\{ y \in Y_i \mid p_i(\gamma^{c_i} \cdot y)
= p_i(y) \big \}\Big) = 0 .\]
\end{enumerate}
Then there is an uncountable set $J \subset I$ such that for every
$i, j \in J$, there are non-null $c_i|_{\F_2}$-invariant and
$c_j|_{\F_2}$-invariant subsets $Y_i'$ and $Y_j'$ of $Y_i$ and
$Y_j$, respectively, so that $c_i|_{\F_2}|Y_i'$ is conjugate to
$c_j|_{\F_2}|Y_j'$.
\end{lemma}

We may obtain the following generalization by changing the
requirement that $\F_2$ is a subgroup of $\Gamma$ to a requirement
$E^{d_i}_{\F_2} \subset E^{c_i}_{\Gamma}$ where $d_i$, $c_i$ are
actions of $\F_2$ and $\Gamma$, respectively. However, only
cosmetic alteration needs to be made to Ioana's original proof.

\begin{lemma} \label{prop:F2-conj} Let $\{c_i \}_{i \in I}$ be an uncountable
set of orbit equivalent ergodic, free, measure preserving actions
$\Gamma \actson^{c_i}(Y_i, m_i)$ such that for each $i$, there is a
free measure preserving action $\F_2 \actson^{d_i} (Y_i, m_i)$ so
that the following conditions hold:
\begin{enumerate}
\item \label{cond:containment} $E^{d_i}_{\F_2} \subset
E^{c_i}_{\Gamma};$

\item \label{cond:quotient} $a$ is a quotient of $d_i$ with
quotient map $p_i: Y_i \rightarrow \T^2;$

\item \label{cond:noclasscollapse} $\forall i \in I \quad \forall
\gamma \in \Gamma \setminus \{ e \}$,
\[m\Big( \big\{ y \in Y_i \mid p_i(\gamma^{c_i} \cdot y)
= p_i(y) \big\}\Big) = 0 .\]
\end{enumerate}
Then there is an uncountable set $J \subset I$ such that for every
$i, j \in J$, there are non-null $d_i$-invariant and $d_j$-invariant
subsets $Y_i' \subset Y_i$ and $Y_j' \subset Y_j$, respectively, so
that $d_i|Y_i'$ is conjugate to $d_j|Y_j'$.
\end{lemma}

Let $\{\pi_i\}_{i \in I}$ be a set of continuum many
non-isomorphic, irreducible weakly mixing representations of
$\F_2$. For each such representation, using Theorem E.1 of
\cite{Kechris2005ap}, obtain a Gaussian action $\F_2
\actson^{a_{\pi_i}} (Z_i, \nu_i)$ such that $\pi_i \cong \pi_j
\implies a_{\pi_i} \cong a_{\pi_j}.$ In addition, we will have
$\pi_i \leq \kappa_0^{a_{\pi_i}}$ and $a_{\pi_i}$ will be a weakly
mixing action.

For each $i \in I$, let the actions $\Gamma \actson^{c_i} (Y_i, m)$
and $\F_2 \actson^{d_i} (Y_i, m)$ and the map $p_i \colon Y_i \to
\T^2 \times Z_i$ be obtained from Lemma \ref{l:action}. Also, let $q
\colon \T^2 \times Z_i \to \T^2$ be given by $q(t, z) = t$. Then
$c_i$, $d_i$ satisfy the conditions of Lemma \ref{prop:F2-conj} with
quotient map $q \circ p_i$.

We claim that for each $i \in I$, the set
\[J_i = \{j \in I \mid c_i
\hbox { is orbit equivalent to } c_j \}\] is countable. Otherwise,
by Lemma \ref{prop:F2-conj}, there is an uncountable set $J \subset
J_i$ such that for any $i, j \in J$, there exist non-null
$\F_2$-invariant subsets $Y_i', Y_j'$ of $Y_i, Y_j$, respectively
such that $d_i|Y_i'$ is orbit equivalent to $d_j|Y_j'$.

If we take the Koopman representation of $d_i$ restricted to $Y_i'$,
then $\kappa_0^{d_i|Y_i'} \leq \kappa_0^{d_i}$. From our
construction, we have
    \[\pi_i \leq \kappa_0^{a_{\pi_i}} \leq \kappa_0^{a \times a_{\pi_i}}\]
since $a_{\pi_i} \sqsubseteq a \times a_{\pi_i}$ and
    \[\kappa_0^{a \times a_{\pi_i}} \leq \kappa_0^{d_i|Y_i'}
\]
since $a \times a_{\pi_i} \sqsubseteq d_i|Y_i'$. As a result,
$\pi_i \leq \kappa_0^{d_i|Y_i'}$ and, finally,
    \[ \pi_i \leq \kappa_0^{d_i|Y_i'} \cong \kappa_0^{d_j|Y_j'}
    \leq \kappa_0^{d_j}.\]
However, a separable unitary representation can only have
countably many non-isomorphic irreducible subrepresentations and
since the $\pi_j$'s are pairwise non-equivalent, it must be that
each $J_i$ is countable and this completes the proof.

\bibliographystyle{alpha}
\bibliography{thesisbib}

\def\cprime{$'$}
\begin{thebibliography}{DGRS08}

\bibitem[BG81]{Bezuglyi1981}
S.~I. Bezuglyi and V.~Ya. Golodets.
\newblock Hyperfinite and {{II}$_{1}$} actions for nonamenable groups.
\newblock {\em J. Funct. Anal.}, 40(1):30--44, 1981.

\bibitem[CFW81]{Connes1981}
A.~Connes, J.~Feldman, and B.~Weiss.
\newblock An amenable equivalence relation is generated by a single
  transformation.
\newblock {\em Ergodic Theory Dynamical Systems}, 1(4):431--450 (1982), 1981.

\bibitem[CW80]{Connes1980}
Alain Connes and Benjamin Weiss.
\newblock Property ({T}) and asymptotically invariant sequences.
\newblock {\em Israel J. Math.}, 37(3):209--210, 1980.

\bibitem[DGRS08]{Dooley2008}
A.~H. Dooley, V.~Ya. Golodets, D.J. Rudolph, and S.D. Sinel'shchikov.
\newblock Non-bernoulli systems with completely positive entropy.
\newblock {\em Ergodic Theory and Dynamical Systems}, 28(1):87--124, 2008.

\bibitem[Dye65]{Dye1965}
H.~A. Dye.
\newblock On the ergodic mixing theorem.
\newblock {\em Trans. Amer. Math. Soc.}, 118:123--130, 1965.

\bibitem[Fer06]{Fernos2006}
Talia Fernos.
\newblock {\em Relative Property (T), Linear Groups, and Applications}.
\newblock PhD thesis, University of Illinois, Chicago, 2006.

\bibitem[FM77]{Feldman1977}
Jacob Feldman and Calvin~C. Moore.
\newblock Ergodic equivalence relations, cohomology, and von {N}eumann
  algebras. {I}.
\newblock {\em Trans. Amer. Math. Soc.}, 234(2):289--324, 1977.

\bibitem[Gab05]{Gaboriau2005b}
D.~Gaboriau.
\newblock Examples of groups that are measure equivalent to the free group.
\newblock {\em Ergodic Theory Dynam. Systems}, 25(6):1809--1827, 2005.

\bibitem[GL]{Gaboriau2007p}
Damien Gaboriau and Russell Lyons.
\newblock A measurable-group-theoretic solution to von {N}eumann's problem.
\newblock preprint, 2007.

\bibitem[Gla03]{Glasner2003}
Eli Glasner.
\newblock {\em Ergodic theory via joinings}, volume 101 of {\em Mathematical
  Surveys and Monographs}.
\newblock American Mathematical Society, Providence, RI, 2003.

\bibitem[GP05]{Gaboriau2005}
Damien Gaboriau and Sorin Popa.
\newblock An uncountable family of nonorbit equivalent actions of
  {$\mathbb{F}\sb n$}.
\newblock {\em J. Amer. Math. Soc.}, 18(3):547--559 (electronic), 2005.

\bibitem[Hjo05]{Hjorth2005}
Greg Hjorth.
\newblock A converse to {D}ye's theorem.
\newblock {\em Trans. Amer. Math. Soc.}, 357(8):3083--3103 (electronic), 2005.

\bibitem[Ioa07]{Ioana2006p}
Adrian Ioana.
\newblock Orbit inequivalent actions for groups containing a copy of
  {$\mathbb{F}\sb 2$}.
\newblock preprint, 2007.

\bibitem[Kec]{Kechris2005ap}
Alexander~S. Kechris.
\newblock Global aspects of ergodic group actions and equivalence relations.
\newblock preprint, 2007.

\bibitem[Kec95]{Kechris1995}
Alexander~S. Kechris.
\newblock {\em Classical descriptive set theory}, volume 156 of {\em Graduate
  Texts in Mathematics}.
\newblock Springer-Verlag, New York, 1995.

\bibitem[KM04]{Kechris2004}
Alexander~S. Kechris and Benjamin~D. Miller.
\newblock {\em Topics in orbit equivalence}, volume 1852 of {\em Lecture Notes
  in Mathematics}.
\newblock Springer-Verlag, Berlin, 2004.

\bibitem[MS06]{Monod2006}
Nicolas Monod and Yehuda Shalom.
\newblock Orbit equivalence rigidity and bounded cohomology.
\newblock {\em Ann. of Math. (2)}, 164(3):825--878, 2006.

\bibitem[MvN36]{Murray1936}
F.~J. Murray and J.~von Neumann.
\newblock On rings of operators.
\newblock {\em Ann. of Math. (2)}, 37(1):116--229, 1936.

\bibitem[OW80]{Ornstein1980}
Donald~S. Ornstein and Benjamin Weiss.
\newblock Ergodic theory of amenable group actions. {I}. {T}he {R}ohlin lemma.
\newblock {\em Bull. Amer. Math. Soc. (N.S.)}, 2(1):161--164, 1980.

\bibitem[Pop06]{Popa2006}
Sorin Popa.
\newblock Strong rigidity of {II$_1$} factors arising from malleable actions of
  {$w$}-rigid groups. {I}.
\newblock {\em Invent. Math.}, 165(2):369--408, 2006.

\bibitem[Sch81]{Schmidt1981}
Klaus Schmidt.
\newblock Amenability, {K}azhdan's property ({T}), strong ergodicity and
  invariant means for ergodic group-actions.
\newblock {\em Ergodic Theory Dynamical Systems}, 1(2):223--236, 1981.

\bibitem[T{\"o}r05]{Tornquist2005}
Asger T{\"o}rnquist.
\newblock {\em The {B}orel complexity of orbit equivalence}.
\newblock PhD thesis, University of California, Los Angeles, 2005.

\bibitem[Zim84]{Zimmer1984}
Robert~J. Zimmer.
\newblock {\em Ergodic theory and semisimple groups}, volume~81 of {\em
  Monographs in Mathematics}.
\newblock Birkh\"auser Verlag, Basel, 1984.

\end{thebibliography}

\end{document}